\documentclass{amsart}
\usepackage{amsmath}
\usepackage{amscd}
\usepackage{amssymb}

\newcommand{\Vol}{\operatorname{Vol}}

\newtheorem{theorem}{Theorem}[section]
\newtheorem{theorem/definition}{Theorem/Definition}[section]
\newtheorem{proposition}{Proposition}[section]
\newtheorem{lemma}{Lemma}[section]
\newtheorem{corollary}{Corollary}[section]
\theoremstyle{remark}
\newtheorem{remark}{Remark}[section]

\theoremstyle{definition}

\begin{document}
\title
{On complete gradient shrinking Ricci solitons}
\author{Huai-Dong Cao and Detang Zhou}
\address{Department of Mathematics\\ Lehigh University\\
Bethlehem, PA 18015\\} \email{huc2@lehigh.edu}
\address{Instituto de Matematica\\ Universidade Federal Fluminense\\
Niter\'oi, RJ 24020\\Brazil} \email{zhou@impa.br}

\begin{abstract} In this paper we derive a precise estimate on the growth
of potential functions of complete noncompact shrinking solitons.
Based on this, we prove that a complete noncompact gradient
shrinking Ricci soliton has at most Euclidean volume growth.
The latter result can be viewed as an analog of the well-known
theorem of Bishop that a complete noncompact Riemannian manifold
with nonnegative Ricci curvature has at most Euclidean volume
growth.

\end{abstract}
\maketitle
\date{}


\footnotetext[1]{The first author was partially supported by NSF grants DMS-0354621 and DMS-0506084;
the second author was partially supported by CNPq and FAPERJ,
Brazil.}

\section{The results}

A complete Riemannian metric $g_{ij}$ on a smooth manifold $M^n$
is called a {\it gradient shrinking Ricci soliton} if there exists
a smooth function $f$ on $M^n$ such that the Ricci tensor $R_{ij}$
of the metric $g_{ij}$ is given by
$$R_{ij}+\nabla_i\nabla_jf=\rho g_{ij}$$
for some positive constant $\rho>0$. The function $f$ is called a
{\it potential function} of the shrinking soliton. Note that by
scaling $g_{ij}$ one can normalize $\rho=\frac{1}{2}$ so that
$$R_{ij}+\nabla_i\nabla_jf=\frac{1}{2} g_{ij}.\eqno(1.1)$$

Gradient shrinking Ricci solitons play an important role in
Hamilton's Ricci flow as they correspond to self-similar
solutions, and often arise as Type I singularity models. In this
paper, we investigate the asymptotic behavior of potential
functions and volume growth rates of complete noncompact gradient
shrinking solitons. Our main results are:

\begin{theorem} Let $(M^n, g_{ij}, f)$ be a complete noncompact
gradient shrinking Ricci soliton satisfying (1.1). Then, the
potential function $f$ satisfies the estimates
$$\frac{1}{4} (r(x)-c_1)^2\leq f(x)\leq \frac{1}{4} (r(x)+c_2)^2.$$
Here $r(x)=d(x_0, x)$ is the distance function from some fixed
point $x_0\in M$, $c_1$ and $c_2$ are positive constants depending
only on $n$ and the geometry of $g_{ij}$ on the unit ball
$B_{x_0}(1)$.
\end{theorem}

\begin{remark} In view of the Gaussian shrinker, namely the flat
Euclidean space $(\mathbb R^n, g_0)$ with the potential function
$|x|^2/4$, the leading term $\frac{1}{4}r^{2}(x)$ for the lower
and upper bounds on $f$ in Theorem 1.1 is optimal. We also point
out that it has been known, by the work of Ni-Wallach \cite{NW1}
and Cao-Chen-Zhu \cite{CCZ08}, that {\sl any 3-dimensional
complete noncompact non-flat shrinking gradient soliton is
necessarily the round cylinder $\mathbb{S}^2\times \mathbb{R}$ or
one of its $\mathbb Z_2$ quotients}.
\end{remark}

\begin{remark} When the Ricci curvature of $(M^n, g_{ij}, f)$ is assumed
to be bounded, Theorem 1.1 was shown by Perelman \cite{P2}. Also,
under the assumption of $Rc\ge 0$, a lower estimate of the
form $$ f(x)\ge \frac{1}{8} r^2(x)-c'_1$$ was shown by Ni
\cite{Ni}. Moreover, the upper bound in Theorem 1.1 was
essentially observed in \cite{CCZ08}, while a rough quadratic
lower bound, as pointed out by Carrillo-Ni \cite{CN},  could follow from the
argument of Fang-Man-Zhang in \cite{FMZ}.
\end{remark}

\begin{theorem} Let $(M^n, g_{ij}, f)$ be a complete noncompact
gradient shrinking Ricci soliton.  Then,
there exists some positive constant $C_1>0$ such that
$$ \Vol(B_{x_0}(r))\leq C_1 r^{n}$$
for $r>0$ sufficiently large.
\end{theorem}

\begin{remark} In an earlier version of the paper, we had an extra
assumption
$$R(x)\le \alpha r^2(x)+ A(r(x)+1), \eqno(1.2)$$ with
$0\le\alpha <\frac{1}{4}$ and $A>0$, on the scalar curvature $R$.
However, as observed by Ovidiu Munteanu, assumption (1.2) actually
is not needed in our proof because there holds estimate (3.4) on
the average of the scalar curvature in general. Note that, as
stated in Lemma 2.3, $R(x)\le \frac{1}{4} (r(x)+c)^2$ holds for
any complete noncompact gradient shrinking soliton. It remains
interesting to find out whether $R$ is bounded from above by a
constant.
\end{remark}

\begin{remark}
Feldman-Ilmanen-Knopf \cite{FIK} constructed a complete noncompact
gradient K\"ahler shrinker on the tautological line bundle
$\mathcal{O}(-1)$ of the complex projective space $\mathbb
CP^{n-1}$ ($n\ge 2$) which has Euclidean volume growth, quadratic
curvature decay, and with Ricci curvature changing signs.
This example shows that the volume growth rate in Theorem 1.2 is
optimal. Note that Carrillo-Ni \cite {CN} showed that {\sl any
non-flat gradient shrinking soliton with nonnegative Ricci
curvature $Rc\ge 0$ must have zero asymptotic volume ratio}, i.e.,
$\lim_{r\to\infty} \Vol(B_{x_0}(r))/r^n=0.$
\end{remark}

Combining Theorem 1.1 and Theorem 1.2, we also have the following
consequence, which was obtained previously in \cite{FM} and
\cite{WW} respectively.

\begin{corollary} Let $(M^n, g_{ij}, f)$ be a complete noncompact
gradient shrinking Ricci soliton.  Then we have $$\int_M |u|
e^{-f}dV< +\infty$$ for any function $u$ on $M$ with $|u(x)|\le
Ae^{\alpha r^2(x)}$, $0\le\alpha <\frac{1}{4}$ and $A>0$. In
particular, the weighted volume of $M$ is finite,
$$\int_M e^{-f}dV < +\infty.$$
\end{corollary}

\noindent {\bf Acknowledgements.} We are grateful to Ovidiu Munteanu for
pointing out to us the estimate (3.4) on the average scalar curvature. The second
author wishes to thank Jiayu Li for some discussions.

\section{Asymptotic behavior of the potential function}

In this section, we investigate the asymptotic behavior of the
potential function of an arbitrary complete noncompact gradient
shrinking Ricci solitons and prove Theorem 1.1.

\smallskip
First of all, we need a few useful facts about complete gradient
shrinking solitons. The first basic result is due to Hamilton (cf.
Theorem 20.1 in \cite{Ha95F}).

\begin{lemma} Let $(M^n, g_{ij}, f)$
be a complete gradient shrinking soliton satisfying (1.1). Then we
have
$$\nabla_iR=2R_{ij}\nabla_jf, $$ and
$$R+|\nabla f|^2-f=C_0 $$ for some constant $C_0$. Here $R$
denotes the scalar curvature of $g_{ij}$.
\end{lemma}

 As a consequence, by adding the constant $C_0$ to $f$, we can assume
$$R+|\nabla f|^2-f=0. \eqno(2.1)$$
From now on we will make this normalization on $f$ throughout the
paper.

\smallskip

We will also need the following useful result, which is a special
case of a more general result on complete ancient solutions due to
B.-L. Chen \cite{BChen} (cf. Proposition 5.5 in \cite{Cao08}).

\begin{lemma} Let $(M^n, g_{ij}, f)$ be a complete shrinking Ricci
soliton. Then $g_{ij}$ has nonnegative scalar curvature $R\ge 0$.
\end{lemma}

As an immediate consequence of (2.1) and Lemma 2.2, one gets the
following result, which was essentially observed by Cao-Chen-Zhu
\cite{CCZ08} (cf. p.78-79 in \cite{CCZ08}).

\begin{lemma} Let $(M^n, g_{ij}, f)$ be a complete shrinking Ricci
soliton satisfying (1.1) and (2.1). Then,
$$f(x)\leq \frac{1}{4}(r(x)+2\sqrt{f(x_0)})^2, \eqno(2.2)$$
$$|\nabla f|(x)\leq \frac{1}{2}r(x)+\sqrt{f(x_0)},
\eqno(2.3)$$ and
$$R(x)\leq \frac{1}{4}(r(x)+2\sqrt{f(x_0)})^2. \eqno(2.4)$$ Here $r(x)=d(x_0, x)$ is the
distance function from some fixed point $x_0\in M$.
\end{lemma}

\begin{proof} By Lemma 2.2 and (2.1),
$$0\leq |\nabla f|^2\leq f, \quad \mbox{or} \quad |\nabla \sqrt{f}|\leq \frac{1}{2} \eqno(2.5)$$ whenever $f>0$.
Thus $\sqrt{f}$ is an Lipschitz function and
$$|\sqrt{f(x)}-\sqrt{f(x_0)}|\leq \frac{1}{2} r(x).$$ Hence
$$\sqrt{f(x)}\leq \frac{1}{2} r(x) + \sqrt{f(x_0)},$$ or
$$f(x)\leq \frac{1}{4} (r(x) +2\sqrt{f(x_0)})^{2}. $$ This proves
(2.2), from which (2.3) and (2.4) follow immediately.
\end{proof}

Now (2.2) provides the upper estimate on $f$ in Theorem 1.1.
However, proving the lower estimate turns out to be more subtle.

\begin{proposition}
Let $(M^n, g_{ij}, f)$ be a complete noncompact gradient shrinking
Ricci soliton satisfying the normalization conditions (1.1) and
(2.1). Then, $f$ satisfies the estimate
$$f(x) \geq \frac{1}{4} (r(x)-c_1)^2,$$
where $c_1$ is a positive constant depending only on $n$ and the
geometry of $g_{ij}$ on the unit ball $B_{x_0}(1)$.
\end{proposition}

\begin{proof}
Consider any minimizing normal geodesic $\gamma(s)$, $0\leq s \leq
s_0$ for some arbitrary large $s_0>0$, starting from
$x_0=\gamma(0)$. Denote by $X(s)=\dot\gamma(s)$ the unit tangent
vector along $\gamma$. Then, by the second variation of arc
length, we have

$$\int_{0}^{s_0} \phi^{2}Rc (X, X) ds \leq (n-1)
\int_{0}^{s_0}|\dot\phi(s)|^2 ds \eqno(2.6)$$
for every nonnegative function $\phi(s)$ defined on the interval
$[0, s_0]$. Now, following Hamilton \cite{Ha95F}, we choose $\phi
(s)$ by
$$ {\phi (s) =\left\{
       \begin{array}{lll}
  s, \ \ \qquad  s\in [0, 1],\\[4mm]
  1, \ \ \qquad  s\in [1, s_0-1], \\[4mm]
  s_0-s, \ \ s\in [s_0-1, s_0].
       \end{array}
    \right.}$$
Then
\begin{align*} \int_{0}^{s_0} Rc (X, X) ds
& = \int_{0}^{s_0} \phi^{2}Rc (X, X)ds + \int_{0}^{s_0}
(1-\phi^{2})Rc (X,
X)ds\\
& \leq (n-1) \int_{0}^{s_0}|\dot\phi(s)|^2 ds + \int_{0}^{s_0}
(1-\phi^{2})Rc (X,
X)ds\\
&\leq 2(n-1) + \max_{B_{x_0}(1)} |Rc| + \max_{B_{\gamma(s_0)}(1)}
|Rc|.
\end{align*}

On the other hand, by (1.1), we have
$$
\nabla_X\dot{f}= \nabla_X\nabla_Xf=\frac{1}{2}- Rc(X,X).
\eqno(2.7)
$$
Integrating (2.7) along $\gamma$ from $0$ to $s_0$, we get
\begin{align*}
\dot{f}(\gamma(s_0))-\dot{f}(\gamma(0))&
=\frac{1}{2}s_{0}-\int_0^{s_0}
Rc(X,X)ds\\
& \geq \frac{s_0}{2} - 2(n-1) - \max_{B_{x_0}(1)} |Rc| -
\max_{B_{\gamma(s_0)}(1)} |Rc|.
\end{align*}

In case $g_{ij}$ has bounded Ricci curvature $|Rc|\le C$ for some
constant $C>0$, then it would follow that
$$
\dot{f}(\gamma(s_0))\geq \frac{1}{2}s_{0}
-\dot{f}(\gamma(0))-2(n-1)-2C=\frac{1}{2}(s_{0}-c), \eqno(2.8)$$
 and that
$$ f(\gamma(s_0)) \geq \frac{1}{4}(s_{0}-c)^2-f(x_0)-\frac{c^2}{4},$$
proving what we wanted.\footnote[2]{Indeed, the above argument was
essentially sketched by Perelman (see, p.3 of \cite{P2}), and a
detailed argument was presented in \cite{CZ05} (p.385-386).}

However, since we do not assume any curvature bound in Theorem
1.1, we have to modify the above argument.

First of all, by integrating (2.7) along $\gamma$ from $s=1$ to $s=s_{0}-1$ instead and using (2.6) as before, we have
\begin{align*}
{\dot{f}}(\gamma(s_{0}-1))-{\dot{f}}(\gamma(1))&
=\int_{1}^{s_0-1}\nabla_X\dot{f}(\gamma(s)) ds\\
&=\frac{1}{2}(s_{0}-2)-\int_{1}^{s_0-1}
Rc(X,X) ds\\
&=\frac{1}{2}(s_{0}-2)-\int_{1}^{s_0-1}
\phi^{2}(s)Rc(X,X) ds\\
& \geq \frac{s_0}{2} - 2n+1 - \max_{B_{x_0}(1)} |Rc| +
\int_{s_0-1}^{s_0} \phi^{2}Rc (X, X) ds.
\end{align*}
Next, using (2.7)  and integration by parts  one more time as in \cite{FMZ},  we
obtain

\begin{align*}\int_{s_0-1}^{s_0} \phi^{2}Rc (X, X) ds &= \frac{1}{2}\int_{s_0-1}^{s_0}
\phi^{2}(s) ds-\int_{s_0-1}^{s_0}
\phi^{2}(s)\nabla_X{\dot{f}}(\gamma(s)) ds\\
&=\frac{1}{6}+{\dot{f}}(\gamma(s_{0}-1))-2\int_{s_0-1}^{s_0}
\phi(s){\dot{f}}(\gamma(s)) ds.
\end{align*}
Therefore,
$$2\int_{s_0-1}^{s_0}
\phi(s){\dot{f}}(\gamma(s)) ds \geq \frac{s_0}{2} - 2n+\frac{7}{6}
 - \max_{B_{x_0}(1)} |Rc| +{\dot{f}}(\gamma(1)).\eqno(2.9)$$
Furthermore, by (2.5) we have
$$|{\dot{f}}(\gamma(s))|\leq \sqrt{f(\gamma(s))},$$ and
$$|\sqrt{f(\gamma(s))}-\sqrt{f(\gamma(s_0))}|\leq \frac{1}{2}
(s_0-s)\leq \frac{1}{2},$$ whenever $s_0-1\leq s \leq s_0$. Thus,
$$\max_{s_0-1\le s\le s_0}|{\dot{f}}(\gamma(s))|\leq
\sqrt{f(\gamma(s_0))}+\frac{1}{2}. \eqno(2.10)$$ Combining (2.9)
and (2.10), and noting $2\int_{s_0-1}^{s_0} \phi(s) ds=1$, we
conclude that
$$\sqrt{f(\gamma(s_0))}\geq \frac{1}{2}(s_0-c_1) $$
for some constant $c_1$ depending only on $n$ and the geometry of
$g_{ij}$ on the unit ball $B_{x_0}(1)$. This completes the proof
of Proposition 2.1 and Theorem 1.1.
\end{proof}

\section{Volume growth of complete gradient shrinking solitons}

In this section, we examine the volume growth of geodesic balls of
complete noncompact gradient shrinking Ricci solitons.

\medskip

Let us define $$\rho(x)=2\sqrt{f(x)}.$$ Then, by Theorem 1.1, we
have
$$r(x)-c\leq \rho(x)\leq r(x)+c \eqno(3.1)$$ with
$c=\max \{c_1, c_2\}>0$. Also, we have
$$\nabla \rho=\frac{\nabla f}{\sqrt{f}} \quad \mbox{and}
\quad |\nabla \rho|=\frac{|\nabla f|}{\sqrt{f}}\le 1.\eqno(3.2)$$

Denote by $$ D(r)=\{x\in M: \rho(x)<r\} \quad \mbox{and} \quad
V(r)=\int_{D(r)}dV. $$ Then, by the Co-Area formula (cf.
\cite{SY}), we have,
$$V(r)=\int_{0}^{r}ds \int_{\partial D(s)}
\frac{1}{|\nabla\rho|} dA. $$ Hence,
$$V'(r)=\int_{\partial D(r)}
\frac{1}{|\nabla\rho|}dA=\frac {r}{2} \int_{\partial D(r)}
\frac{1}{|\nabla {f}|}dA. \eqno(3.3)$$ Here we have used (3.2) in
deriving the last identity in (3.3).
\begin{lemma}
$$n V(r)- r V'(r) = 2 \int_{D(r)}{R}dV - 2
\int_{\partial D(r)}\frac{R}{|\nabla f|}dV.$$
\end{lemma}

\begin{proof}
Taking the trace in (1.1), we have
$$R+\Delta f=\frac{n}{2}.$$ Thus,
\begin{align*}
 n V(r)-2 \int_{D(r)}{R}dV & =
2\int_{D(r)} \Delta f dV\\
& = 2 \int_{\partial D(r)} \nabla f\cdot\frac{\nabla
\rho}{|\nabla\rho|}\\
& =2 \int_{\partial D(r)} |\nabla f|dV\\
& = 2 \int_{\partial D(r)}\frac{f-R} {|\nabla f|}dV\\
 & =r V'(r)- 2\int_{\partial D(r)}\frac{R}{|\nabla f|}dV.
\end{align*}

\end{proof}

\begin{remark} As pointed out to us by Ovidiu Munteanu, we have
also actually shown that $$ \int_{D(r)}{R}dV \le \frac{n}{2}
V(r).\eqno(3.4)$$ Namely, the average scalar curvature over $D(r)$
is bounded by $n/2$.

\end{remark}

Now we are ready to prove {\bf Theorem 1.2}.

\medskip
\begin{proof} Let $(M^n, g_{ij}, f)$ be a complete noncompact gradient
shrinking Ricci soliton.  Denote
$$ \chi (r)=\int_{D(r)}RdV. $$ By the Co-Area formula, we have
$$\chi(r) = \int_{0}^r ds\int_{\partial D(s)} \frac{R}{|\nabla \rho|}dA
        = \frac12\int_{0}^rsds\int_{\partial D(s)} \frac{R}{|\nabla
{f}|}dA. $$
Hence
$$
  \chi'(r)=\frac  r2\int_{\partial D(r)} \frac{R}{|\nabla
{f}|}dA. $$

Therefore, Lemma 3.1 can be rewritten as
$$n V(r)- r V'(r) = 2\chi(r) - \frac{4}{r}\chi'(r). \eqno(3.5)$$
This implies that
\begin{align*}
(r^{-n}V(r))' & =4r^{-n-2}e^{\frac{r^2}4}(e^{-\frac{r^2}4}\chi(r))'\\
& =4r^{-n-2}\chi'(r)-2r^{-n-1}\chi(r).
\end{align*}
Integrating the above equation from $r_0$ to $r$, we get
\begin{align*}
r^{-n}V(r)-r_0^{-n}V(r_0) &=\left.4r^{-n-2}\chi(r)\right|^{r}_{r_0}+4(n+2)\int_{r_0}^{r} r^{-n-3}\chi(r)dr\\
         &\qquad -2\int_{r_0}^{r} r^{-n-1}\chi(r)dr\\
         &=4r^{-n-2}\chi(r)-4r_0^{-n-2}\chi(r_0)\\
         &\qquad + 2\int_{r_0}^{r} r^{-n-3}\chi(r)(2(n+2)-r^2)dr.
\end{align*}
Since $\chi(r)$ is positive and increasing in $r$ we have, for
$r_0= \sqrt{2(n+2)}$,
\begin{align*}
\int_{r_0}^{r} r^{-n-3}\chi(r)(2(n+2)-r^2)dr & \le
\chi(r_0)\int_{r_0}^{r} r^{-n-3}(2(n+2)-r^2)dr\\
&=\chi(r_0)\left.(-2r^{-n-2}+\frac{1}{n} r^{-n})\right|^{r}_{r_0}.
\end{align*}
Thus,
$$r^{-n}V(r)-r_0^{-n}V(r_0)\le 4r^{-n-2}(\chi(r)-\chi(r_0))+
\frac{2}{n} \chi(r_0) (r^{-n}-r_{0}^{-n}).$$ It follows that, for
$r\ge r_0= \sqrt{2(n+2)}$,
$$V(r)\le  (r_0^{-n}V(r_0))r^{n}+
4r^{-2}\chi(r).\eqno(3.6)$$ On the other hand, by (3.4) we have
$$4r^{-2}\chi(r)\le 2n r^{-2} V(r)\le \frac{1}{2} V(r), \eqno(3.7)$$
for $r$ sufficiently large.

Plugging (3.7) into (3.6), we obtain
$$V(r)\le 2r_0^{-n}V(r_0)r^{n}.$$ Therefore, by (3.1),
$$\Vol(B_{x_0}(r))\leq V(r+c)\leq  V(r_0)r^{n}$$  for $r$ large enough.
This finishes the proof of
Theorem 1.2.

\end{proof}

We conclude with the following volume lower estimate.

\begin{proposition} Let $(M^n, g_{ij}, f)$ be a complete noncompact
gradient shrinking Ricci soliton. Suppose the average scalar
curvature satisfies the upper bound
$$ \frac{1}{V(r)}\int_{D(r)}{R}dV \le \delta, \eqno(3.8)$$
for some positive constant $\delta < n/2$ and all sufficiently
large $r$. Then, there exists some positive constant $C_2>0$ such
that
$$\Vol(B_{x_0}(r))\geq C_2 r^{n-2\delta}$$
for $r$ sufficiently large.
\end{proposition}

\begin{proof} Combining the assumption (3.8) with Lemma 3.1 and Lemma 2.2, it follows that
$$(n-2\delta) V(r)\leq rV'(r). \eqno(3.9)$$
Thus,
$$\int_{1}^{r}
\frac{V'(s)}{V(s)} ds \geq \int_{1}^r\frac{n-2\delta}{s} ds.$$
Consequently,
$$V(r)\geq V(1) r^{n-2\delta}.$$ Therefore, in view of (3.1),
$$\Vol(B_{x_0}(r))\geq V(r-c)\geq 2^{-n}V(1)r^{n-2\delta}$$ for $r$ sufficiently large.

\end{proof}

\begin{remark}
X.-P. Zhu and the first author (see Theorem 3.1 in \cite{Cao09})
have shown that a complete noncompact gradient shrinking soliton,
without any curvature assumption, must have infinite volume. Their proof is, however, more sophisticated, relying on a
logarithmic inequality of Carrillo-Ni \cite{CN} and the Perelman
type noncollapsing argument for complete gradient shrinking
solitons.
\end{remark}

\end{document}